\documentclass[a4paper]{amsart}
\usepackage[T1]{fontenc}
\usepackage[utf8]{inputenc}
\usepackage[active]{srcltx}
\usepackage{amsthm}
\usepackage{amstext}
\usepackage{amssymb}
\usepackage{setspace}
\usepackage{xargs}[2008/03/08]

\makeatletter

\newcommand{\noun}[1]{\textsc{#1}}

\theoremstyle{plain}
\newtheorem{thm}{\protect\theoremname}
  \theoremstyle{remark}
  \newtheorem*{rem*}{\protect\remarkname}

\makeatother

  \providecommand{\remarkname}{Remark}
\providecommand{\theoremname}{Theorem}

\begin{document}

\title{Poly-Cauchy and Peters mixed-type polynomials}

\author{Dae San Kim, Taekyun Kim}
\begin{abstract}
The Peters polynomials are a generalization of Boole polynomials. In
this paper, we consider Peters and poly-Cauchy mixed type polynomials
and investigate the properties of those polynomials which are derived
from umbral calculus. Finally, we give various identities of those
polynomials associated with special polynomials.

\newcommandx\li[1][usedefault, addprefix=\global, 1=]{\textnormal{Li}_{#1}}

\newcommandx\lif[1][usedefault, addprefix=\global, 1=]{\textnormal{Lif}_{#1}}

\global\long\def\ch{\textnormal{Ch}_{n}\left(x\right)}

\end{abstract}
\maketitle

\section{Introduction}

The Peters polynomials are defined by the generating function to be
\begin{equation}
\sum_{n=0}^{\infty}S_{n}\left(x;\lambda,\mu\right)\frac{t^{n}}{n!}=\left(1+\left(1+t\right)^{\lambda}\right)^{-\mu}\left(1+t\right)^{x},\quad\left(\textrm{see [14]}\right).\label{eq:1}
\end{equation}

The first few of them are given by
\[
S_{0}\left(x;\lambda,\mu\right)=2^{-\mu},\quad S_{1}\left(x;\lambda,\mu\right)=2^{-\left(\mu+1\right)}\left(2x-\lambda\mu\right),\,\cdots.
\]

If $\mu=1$, then $S_{n}\left(x;\lambda\right)=S_{n}\left(x:\lambda,1\right)$
are called Boole polynomials.

In particular, for $\mu=1$, $\lambda=1$ , $S_{n}\left(x;1,1\right)=\ch$
are Changhee polynomials which are defined by
\[
\sum_{n=0}^{\infty}\ch\frac{t^{n}}{n!}=\frac{1}{t+2}\left(1+t\right)^{x},\quad\left(\textrm{see \ensuremath{\left[8\right]}}\right).
\]

The generating functions for the poly-Cauchy polynomials of the first
kind $C_{n}^{\left(k\right)}\left(x\right)$ and of the second kind
$\hat{C}_{n}^{\left(k\right)}\left(x\right)$ are, respectively, given
by
\begin{equation}
\lif[k]\left(\log\left(1+t\right)\right)\left(1+t\right)^{-x}=\sum_{n=0}^{\infty}C_{n}^{\left(k\right)}\left(x\right)\frac{t^{n}}{n!},\label{eq:2}
\end{equation}
and
\begin{equation}
\lif[k]\left(-\log\left(1+t\right)\right)\left(1+t\right)^{x}=\sum_{n=0}^{\infty}\hat{C}_{n}^{\left(k\right)}\left(x\right)\frac{t^{n}}{n!},\label{eq:3}
\end{equation}
where $\lif[k]\left(t\right)={\displaystyle \sum_{n=0}^{\infty}\frac{t^{n}}{n!\left(n+1\right)^{k}},\,\left(k\in\mathbb{Z}\right),\,\left(\textrm{see }[10,11]\right).}$

In this paper, we consider the poly-Cauchy of the first kind and Peters
(respectively the poly-Cauchy of the second kind and Peters) mixed-type
polynomials as follows:
\begin{equation}
\left(1+\left(1+t\right)^{\lambda}\right)^{-\mu}\lif[k]\left(\log\left(1+t\right)\right)\left(1+t\right)^{-x}=\sum_{n=0}^{\infty}CP_{n}^{\left(k\right)}\left(x;\lambda,\mu\right)\frac{t^{n}}{n!},\label{eq:4}
\end{equation}
 and
\begin{equation}
\left(1+\left(1+t\right)^{\lambda}\right)^{-\mu}\lif[k]\left(-\log\left(1+t\right)\right)\left(1+t\right)^{x}=\sum_{n=0}^{\infty}\hat{C}P_{n}^{\left(k\right)}\left(x;\lambda,\mu\right)\frac{t^{n}}{n!}.\label{eq:5}
\end{equation}

\begin{singlespace}
For $\alpha\in\mathbb{Z}_{\ge0}$, the Bernoulli polynomials
of order $\alpha$ are defined by the generating function to be
\begin{equation}
\left(\frac{t}{e^{t}-1}\right)^{\alpha}e^{xt}=\sum_{n=0}^{\infty}B_{n}^{\left(\alpha\right)}\left(x\right)\frac{t^{n}}{n!},\quad\left(\textrm{see [1-8]}\right).\label{eq:6}
\end{equation}

\end{singlespace}

As is well known, the Frobenius-Euler polynomials of order $\alpha$
are given by
\begin{equation}
\left(\frac{1-\lambda}{e^{t}-\lambda}\right)^{\alpha}e^{xt}=\sum_{n=0}^{\infty}H_{n}^{\left(\alpha\right)}\left(x|\lambda\right)\frac{t^{n}}{n!},\quad\left(\textrm{see [1-13]}\right),\label{eq:7}
\end{equation}
 where $\lambda\in\mathbb{C}$ with $\lambda\ne1$ and $\alpha\in\mathbb{Z}_{\ge0}$.

When $x=0$, $CP_{n}^{\left(k\right)}\left(0;\lambda,\mu\right)$(or
$\hat{C}P_{n}^{\left(k\right)}\left(0;\lambda,\mu\right)$) are called
the poly-Cauchy of the first kind and Peters (or the poly-Cauchy of
the second kind and Peters) mixed-type numbers.

The higher-order Cauchy polynomials of the first kind are defined
by the generating function to be
\begin{equation}
\left(\frac{t}{\log\left(1+t\right)}\right)^{\alpha}\left(1+t\right)^{-x}=\sum_{n=0}^{\infty}\mathbb{C}_{n}^{\left(\alpha\right)}\left(x\right)\frac{t^{n}}{n!},\quad\left(\alpha\in\mathbb{Z}_{\ge0}\right),\label{eq:8}
\end{equation}
 and the higher-order Cauchy polynomials of the second kind are given
by
\begin{equation}
\left(\frac{t}{\left(1+t\right)\log\left(1+t\right)}\right)^{\alpha}\left(1+t\right)^{x}=\sum_{n=0}^{\infty}\mathbb{\hat{C}}_{n}^{\left(\alpha\right)}\left(x\right)\frac{t^{n}}{n!},\quad\left(\alpha\in\mathbb{Z}_{\ge0}\right).\label{eq:9}
\end{equation}

The Stirling number of the first kind is given by
\begin{equation}
\left(x\right)_{n}=x\left(x-1\right)\cdots\left(x-n+1\right)=\sum_{l=0}^{n}S_{1}\left(n,l\right)x^{l}.\label{eq:10}
\end{equation}

Thus, by (\ref{eq:10}), we get
\begin{equation}
\left(\log\left(1+t\right)\right)^{m}=m!\sum_{l=m}^{\infty}S_{1}\left(l,m\right)\frac{t^{l}}{l!},\quad m\in\mathbb{Z}_{\ge0},\;\left(\textrm{see \ensuremath{\left[14\right]}}\right).\label{eq:11}
\end{equation}

It is easy to show that
\begin{equation}
x^{\left(n\right)}=x\left(x+1\right)\cdots\left(x+n-1\right)=\left(-1\right)^{n}\left(-x\right)_{n}=\sum_{l=0}^{n}S_{1}\left(n,l\right)\left(-1\right)^{n-l}x^{l}.\label{eq:12}
\end{equation}

Let $\mathbb{C}$ be the complex number field and let $\mathcal{F}$
be the algebra of all formal power series in the variable $t$ over $\mathbb{C}$
as follows:
\begin{equation}
\mathcal{F}=\left\{ \left.f\left(t\right)=\sum_{k=0}^{\infty}a_{k}\frac{t^{k}}{k!}\right|a_{k}\in\mathbb{C}\right\} .\label{eq:13}
\end{equation}

\begin{singlespace}
Let $\mathbb{P}=\mathbb{C}\left[x\right]$ and let $\mathbb{P}^{*}$
be the vector space of all linear functionals on $\mathbb{P}$. $\left\langle \left.L\right|p\left(x\right)\right\rangle $
denotes the action of linear functional $L$ on the polynomial $p\left(x\right)$,
and we recall that the vector space operations on $\mathbb{P}^{*}$
are defined by $\left\langle \left.L+M\right|p\left(x\right)\right\rangle =\left\langle \left.L\right|p\left(x\right)\right\rangle +\left\langle \left.M\right|p\left(x\right)\right\rangle $,
$\left\langle \left.cL\right|p\left(x\right)\right\rangle =c\left\langle \left.L\right|p\left(x\right)\right\rangle $,
where $c$ is complex constant in $\mathbb{C}$.
\end{singlespace}

For $f\left(t\right)\in\mathcal{F}$, let us define the linear functional
on $\mathbb{P}$ by setting
\begin{equation}
\left\langle \left.f\left(t\right)\right|x^{n}\right\rangle =a_{n},\quad\left(n\ge0\right).\label{eq:14}
\end{equation}

\begin{singlespace}
Then, by (\ref{eq:13}) and (\ref{eq:14}), we get
\begin{equation}
\left\langle \left.t^{k}\right|x^{n}\right\rangle =n!\delta_{n,k},\quad\left(n,k\ge0\right),\quad\left(\textrm{see \ensuremath{\left[14,15\right]}}\right),\label{eq:15}
\end{equation}
where $\delta_{n,k}$ is the Kronecker's symbol.

Let $f_{L}\left(t\right)={\displaystyle \sum_{k=0}^{\infty}\frac{\left\langle \left.L\right|x^{k}\right\rangle }{k!}t^{k}.}$
Then, by (\ref{eq:14}), we see that $\left\langle \left.f_{L}\left(t\right)\right|x^{n}\right\rangle =\left\langle \left.L\right|x^{n}\right\rangle .$
The map $L\longmapsto f_{L}\left(t\right)$ is a vector space isomorphism
from $\mathbb{P}^{*}$ onto $\mathcal{F}$. Henceforth, $\mathcal{F}$
denotes both the algebra of formal power series in $t$ and vector
space of all linear functionals on $\mathbb{P}$, and so an element
$f\left(t\right)$ of $\mathcal{F}$ will be thought of as both a
formal power series and a linear functional. We call $\mathcal{F}$
the umbral algebra and the umbral calculus is the study of umbral
algebra. The order $O\left(f\right)$ of the power series $f\left(t\right)\left(\ne0\right)$
is the smallest integer for which the coefficient of $t^{k}$ does
not vanish. If $O\left(f\left(t\right)\right)=1$, then $f\left(t\right)$
is called a delta series; if $O\left(f\left(t\right)\right)=0$, then
$f\left(t\right)$ is called an invertible series. For $f\left(t\right),g\left(t\right)\in\mathcal{F}$
with $O\left(f\left(t\right)\right)=1$ and $O\left(g\left(t\right)\right)=0$,
there exists a unique sequence $s_{n}\left(x\right)$ of polynomials
such that $\left\langle \left.g\left(t\right)f\left(t\right)^{k}\right|s_{n}(x)\right\rangle =n!\delta_{n,k},\quad\left(n,k\ge0\right)$.

The sequence $s_{n}\left(x\right)$ is called the Sheffer sequence for
$\left(g\left(t\right),f\left(t\right)\right)$ which is denoted by
$s_{n}\left(x\right)\sim\left(g\left(t\right),f\left(t\right)\right).$

For $f\left(t\right),g\left(t\right)\in\mathcal{F}$ and $p\left(x\right)\in\mathbb{P},$we
have
\[
\left\langle \left.f\left(t\right)g\left(t\right)\right|p\left(x\right)\right\rangle =\left\langle \left.f\left(t\right)\right|g\left(t\right)p\left(x\right)\right\rangle =\left\langle \left.g\left(t\right)\right|f\left(t\right)p\left(x\right)\right\rangle ,
\]
and
\begin{equation}
f\left(t\right)=\sum_{k=0}^{\infty}\left\langle \left.f\left(t\right)\right|x^{k}\right\rangle \frac{t^{k}}{k!},\quad p\left(x\right)=\sum_{k=0}^{\infty}\left\langle \left.t^{k}\right|p\left(x\right)\right\rangle \frac{x^{k}}{k!},\label{eq:16}
\end{equation}
By (\ref{eq:16}), we get
\begin{equation}
p^{\left(k\right)}\left(0\right)=\left\langle \left.t^{k}\right|p\left(x\right)\right\rangle =\left\langle \left.1\right|p^{\left(k\right)}\left(x\right)\right\rangle ,\:\left(k\ge0\right).\label{eq:17}
\end{equation}
where $p^{\left(k\right)}\left(0\right)=\left.\frac{d^{k}p\left(x\right)}{dx^{k}}\right|_{x=0}.$
\end{singlespace}

Thus, by (\ref{eq:17}), we have
\begin{equation}
t^{k}p\left(x\right)=p^{\left(k\right)}\left(x\right)=\frac{d^{k}p\left(x\right)}{dx^{k}},\quad\left(\textrm{see \ensuremath{\left[8,10,14\right]}}\right).\label{eq:18}
\end{equation}

Let $s_{n}\left(x\right)\sim\left(g\left(t\right),f\left(t\right)\right)$.
Then the following equations are known in {[}14{]} :
\begin{equation}
\frac{1}{g\left(\overline{f}\left(t\right)\right)}e^{x\overline{f}\left(t\right)}=\sum_{n=0}^{\infty}s_{n}\left(x\right)\frac{t^{n}}{n!},\textrm{ for all }x\in\mathbb{C},\label{eq:19}
\end{equation}
where $\overline{f}\left(t\right)$ is the compositional inverse for
$f\left(t\right)$ with $f\left(\overline{f}\left(t\right)\right)=t$,
\begin{equation}
s_{n}\left(x\right)=\sum_{j=0}^{n}\frac{1}{j!}\left\langle \left.\frac{\left(\overline{f}\left(t\right)\right)^{j}}{g\left(\overline{f}\left(t\right)\right)}\right|x^{n}\right\rangle x^{j},\label{eq:20}
\end{equation}
\begin{equation}
s_{n}\left(x+y\right)=\sum_{j=0}^{n}\dbinom{n}{j}s_{j}\left(x\right)P_{n-j}\left(y\right),\quad\textrm{where }P_{n}\left(x\right)=g\left(t\right)s_{n}\left(x\right),\label{eq:21}
\end{equation}
and
\begin{equation}
s_{n+1}\left(x\right)=\left(x-\frac{g^{\prime}\left(t\right)}{g\left(t\right)}\right)\frac{1}{f^{\prime}\left(t\right)}s_{n}\left(x\right),\quad f\left(t\right)s_{n}\left(x\right)=ns_{n-1}\left(x\right),\quad\left(n\ge0\right),\label{eq:22}
\end{equation}
 and
\begin{equation}
\frac{d}{dx}s_{n}\left(x\right)=\sum_{l=0}^{n-1}\dbinom{n}{l}\left\langle \left.\overline{f}\left(t\right)\right|x^{n-l}\right\rangle s_{l}\left(x\right).\label{eq:23}
\end{equation}

As is well known, the transfer formula for $p_{n}\left(x\right)\sim\left(1,f\left(t\right)\right)$,
$q_{n}\left(x\right)\sim\left(1,g\left(t\right)\right),$ $\left(n\ge1\right)$,
is given by
\begin{equation}
q_{n}\left(x\right)=x\left(\frac{f\left(t\right)}{g\left(t\right)}\right)^{n}x^{-1}p_{n}\left(x\right).\label{eq:24}
\end{equation}

For $s_{n}\left(x\right)\sim\left(g\left(t\right),f\left(t\right)\right),\, r_{n}\left(x\right)\sim\left(h\left(t\right),l\left(t\right)\right),$
let
\begin{equation}
s_{n}\left(x\right)=\sum_{m=0}^{\infty}C_{n,m}r_{n}\left(x\right),\label{eq:25}
\end{equation}
where
\begin{equation}
C_{n,m}=\frac{1}{m!}\left\langle \left.\frac{h\left(\overline{f}\left(t\right)\right)}{g\left(\overline{f}\left(t\right)\right)}\left(l\left(\overline{f}\left(t\right)\right)\right)^{m}\right|x^{n}\right\rangle ,\textrm{ (see [14]).}\label{eq:26}
\end{equation}

It is known that

\begin{equation}
\left\langle \left.f\left(t\right)\right|xp\left(x\right)\right\rangle =\left\langle \left.\partial_{t}f\left(t\right)\right|p\left(x\right)\right\rangle ,\, e^{yt}p\left(x\right)=p\left(x+y\right)\label{eq:27}
\end{equation}
where $f\left(t\right)\in\mathcal{F}$ and $p\left(x\right)\in\mathbb{P}$
(see {[}8,10,14{]}).

In this paper, we consider Peters and poly-Cauchy mixed type polynomials
with umbral calculus viewpoint and investigate the properties of
those polynomials which are derived from umbral calculus. Finally,
we give some interesting identities of those polynomials associated
with speical polynomials.

\begin{singlespace}

\section{poly-Cauchy and Peters mixed-type polynomials}
\end{singlespace}

From (\ref{eq:2}), (\ref{eq:3}), and (\ref{eq:19}), we note that
\begin{equation}
CP_{n}^{\left(k\right)}\left(x;\lambda,\mu\right)\sim\left(\left(1+e^{-\lambda t}\right)^{\mu}\frac{1}{\lif[k]\left(-t\right)},e^{-t}-1\right)\label{eq:28}
\end{equation}
and
\begin{equation}
\hat{C}P_{n}^{\left(k\right)}\left(x;\lambda,\mu\right)\sim\left(\left(1+e^{\lambda t}\right)^{\mu}\frac{1}{\lif[k]\left(-t\right)},e^{t}-1\right).\label{eq:29}
\end{equation}

It is not difficult to show that

\begin{align}
 & \left(1+e^{-\lambda t}\right)^{\mu}\label{eq:30}\\
= & 2^{\mu}\left(1+\frac{1}{2}\sum_{j=1}^{\infty}\frac{\left(-\lambda t\right)^{j}}{j!}\right)^{\mu}\nonumber \\
= & \sum_{i=0}^{\infty}\sum_{j=0}^{\infty}\sum_{j_{1}+\cdots+j_{i}=j}2^{\mu-i}\dbinom{\mu}{i}\dbinom{j+i}{j_{1}+1,\cdots,j_{i}+1}\frac{\left(-\lambda t\right)^{j+i}}{\left(j+i\right)!}\nonumber
\end{align}

and
\begin{align}
 & \left(1+\left(1+t\right)^{\lambda}\right)^{-\mu}\label{eq:31}\\
= & 2^{-\mu}\left(1+\frac{1}{2}\sum_{j=0}^{\infty}\dbinom{\lambda}{j+1}t^{j+1}\right)^{-\mu}\nonumber \\
= & \sum_{i=0}^{\infty}\sum_{j=0}^{\infty}\sum_{j_{1}+\cdots+j_{i}=j}2^{-\left(\mu+i\right)}\dbinom{-\mu}{i}\dbinom{\lambda}{j_{1}+1}\cdots\dbinom{\lambda}{j_{i}+1}t^{j+i}.\nonumber
\end{align}

From (\ref{eq:14}), we have
\begin{align}
 & CP_{n}^{\left(k\right)}\left(y;\lambda,\mu\right)\label{eq:32}\\
= & \left\langle \left.\sum_{l=0}^{\infty}CP_{l}^{\left(k\right)}\left(y;\lambda,\mu\right)\frac{t^{l}}{l!}\right|x^{n}\right\rangle \nonumber \\
= & \left\langle \left.\left(1+\left(1+t\right)^{\lambda}\right)^{-\mu}\lif[k]\left(\log\left(1+t\right)\right)\left(1+t\right)^{-y}\right|x^{n}\right\rangle \nonumber \\
= & \left\langle \left.\left(1+\left(1+t\right)^{\lambda}\right)^{-\mu}\right|\sum_{l=0}^{n}\dbinom{n}{l}C_{l}^{\left(k\right)}\left(y\right)x^{n-l}\right\rangle \nonumber \\
= & \sum_{l=0}^{n}\dbinom{n}{l}C_{l}^{\left(k\right)}\left(y\right)\left\langle \left.\sum_{m=0}^{\infty}S_{m}\left(0;\lambda,\mu\right)\frac{t^{m}}{m!}\right|x^{n-l}\right\rangle \nonumber \\
= & \sum_{l=0}^{n}\dbinom{n}{l}S_{n-l}\left(0;\lambda,\mu\right)C_{l}^{\left(k\right)}\left(y\right).\nonumber
\end{align}

Therefore, by (\ref{eq:32}), we obtain the following theorem.
\begin{thm}
For $n\ge0$, we have
\[
CP_{n}^{\left(k\right)}\left(x;\lambda,\mu\right)=\sum_{l=0}^{n}\dbinom{n}{l}S_{n-l}\left(0;\lambda,\mu\right)C_{l}^{\left(k\right)}\left(x\right).
\]

\end{thm}
$\,$

Alternatively,
\begin{eqnarray}
CP_{n}^{\left(k\right)}\left(y;\lambda,\mu\right) & = & \left\langle \left.\sum_{l=0}^{\infty}CP_{l}^{\left(k\right)}\left(y;\lambda,\mu\right)\frac{t^{l}}{l!}\right|x^{n}\right\rangle \label{eq:33}\\
 & = & \left\langle \left.\lif[k]\left(\log\left(1+t\right)\right)\right|\left(1+\left(1+t\right)^{\lambda}\right)^{-\mu}\left(1+t\right)^{-y}x^{n}\right\rangle \nonumber \\
 & = & \left\langle \left.\lif[k]\left(\log\left(1+t\right)\right)\right|\sum_{l=0}^{n}\dbinom{n}{l}S_{l}\left(-y;\lambda,\mu\right)x^{n-l}\right\rangle \nonumber \\
 & = & \sum_{l=0}^{n}\dbinom{n}{l}S_{l}\left(-y;\lambda,\mu\right)\left\langle \left.\lif[k]\left(\log\left(1+t\right)\right)\right|x^{n-l}\right\rangle \nonumber \\
 & = & \sum_{l=0}^{n}\dbinom{n}{l}S_{l}\left(-y;\lambda,\mu\right)C_{n-l}^{\left(k\right)}\left(0\right).\nonumber
\end{eqnarray}

Therefore, by (\ref{eq:33}), we obtain the following theorem.
\begin{thm}
For $n\ge0$, let $C_{n-l}^{\left(k\right)}\left(0\right)=C_{n-l}^{\left(k\right)}.$
Then we have
\[
CP_{n}^{\left(k\right)}\left(x;\lambda,\mu\right)=\sum_{l=0}^{n}\dbinom{n}{l}C_{n-l}^{\left(k\right)}S_{l}\left(-x;\lambda,\mu\right).
\]
\end{thm}
\begin{rem*}
By the same method, we get
\begin{equation}
\hat{C}P_{n}^{\left(k\right)}\left(x;\lambda,\mu\right)=\sum_{l=0}^{n}\dbinom{n}{l}S_{n-l}\left(0;\lambda,\mu\right)\hat{C}_{l}^{\left(k\right)}\left(x\right),\label{eq:34}
\end{equation}
 and
\begin{equation}
\hat{C}P_{n}^{\left(k\right)}\left(x;\lambda,\mu\right)=\sum_{l=0}^{n}\dbinom{n}{l}\hat{C}_{n-l}^{\left(k\right)}S_{l}\left(x;\lambda,\mu\right).\label{eq:35}
\end{equation}

\end{rem*}
$\,$

From (\ref{eq:20}) and (\ref{eq:28}), we have
\begin{align}
 & CP_{n}^{\left(k\right)}\left(x;\lambda,\mu\right)\label{eq:36}\\
= & \sum_{j=0}^{n}\frac{1}{j!}\left\langle \left.\left(1+\left(1+t\right)^{\lambda}\right)^{-\mu}\lif[k]\left(\log\left(1+t\right)\right)\left(-\log\left(1+t\right)\right)^{j}\right|x^{n}\right\rangle x^{j}\nonumber
\end{align}

From (\ref{eq:31}), we note that
\begin{align}
 & \left\langle \left.\left(1+\left(1+t\right)^{\lambda}\right)^{-\mu}\lif[k]\left(\log\left(1+t\right)\right)\left(-\log\left(1+t\right)\right)^{j}\right|x^{n}\right\rangle \label{eq:37}\\
= & \sum_{m=0}^{n-j}\frac{\left(-1\right)^{j}}{m!\left(m+1\right)^{k}}\sum_{l=0}^{n-j-m}\frac{\left(m+j\right)!}{\left(l+m+j\right)!}S_{1}\left(l+m+j,m+j\right)\nonumber \\
 & \times\left(n\right)_{l+m+j}\left\langle \left.\left(1+\left(1+t\right)^{\lambda}\right)^{-\mu}\right|x^{n-l-m-j}\right\rangle \nonumber \\
= & \sum_{m=0}^{n-j}\frac{\left(-1\right)^{j}}{m!\left(m+1\right)^{k}}\sum_{l=0}^{n-j-m}\frac{\left(m+j\right)!}{\left(l+m+j\right)!}\nonumber \\
 & \times S_{1}\left(l+m+j,m+j\right)\left(n\right)_{l+m+j}\sum_{i=0}^{n-j-m-l}\sum_{r=0}^{\infty}\sum_{r_{1}+\cdots+r_{i}=r}2^{-\left(\mu+i\right)}\nonumber \\
 & \times\dbinom{-\mu}{i}\dbinom{\lambda}{r_{1}+1}\cdots\dbinom{\lambda}{r_{i}+1}\left\langle \left.t^{r+i}\right|x^{n-l-m-j}\right\rangle \nonumber \\
= & 2^{-\mu}n!\sum_{m=0}^{n-j}\sum_{l=0}^{n-j-m}\sum_{i=0}^{n-j-m-l}\sum_{r_{1}+\cdots+r_{i}=n-j-m-l-i}\frac{2^{-i}\left(-1\right)^{j}\left(m+j\right)!}{m!\left(m+1\right)^{k}\left(l+m+j\right)!}\nonumber \\
 & \times\dbinom{-\mu}{i}\dbinom{\lambda}{r_{1}+1}\cdots\dbinom{\lambda}{r_{i}+1}S_{1}\left(l+m+j,m+j\right).\nonumber
\end{align}

Therefore, by (\ref{eq:36}) and (\ref{eq:37}), we obtain the following
theorem.
\begin{thm}
\label{thm:3}For $n\ge0$, we have
\begin{align*}
 & CP_{n}^{\left(k\right)}\left(x;\lambda,\mu\right)\\
= & 2^{-\mu}n!\sum_{j=0}^{n}\frac{\left(-1\right)^{j}}{j!}\left\{ \sum_{m=0}^{n-j}\sum_{l=0}^{n-j-m}\sum_{i=0}^{n-j-m-l}\sum_{r_{1}+\cdots+r_{i}=n-j-m-l-i}\frac{2^{-i}}{m!\left(m+1\right)^{k}}\right.\\
 & \left.\times\frac{\left(m+j\right)!}{\left(l+m+j\right)!}\dbinom{-\mu}{i}\dbinom{\lambda}{r_{1}+1}\cdots\dbinom{\lambda}{r_{i}+1}S_{1}\left(l+m+j,m+j\right)\right\} x^{j}.
\end{align*}
\end{thm}
\begin{rem*}
By the same method as Theorem \ref{thm:3}, we get
\begin{align}
 & \hat{C}P_{n}^{\left(k\right)}\left(x;\lambda,\mu\right)\label{eq:38}\\
= & 2^{-\mu}n!\sum_{j=0}^{n}\frac{1}{j!}\left\{ \sum_{m=0}^{n-j}\sum_{l=0}^{n-j-m}\sum_{i=0}^{n-j-m-l}\sum_{r_{1}+\cdots+r_{i}=n-j-m-l-i}\frac{2^{-i}\left(-1\right)^{m}}{m!\left(m+1\right)^{k}}\right.\nonumber \\
 & \times\frac{\left(m+j\right)!}{\left(l+m+j\right)!}\dbinom{-\mu}{i}\nonumber \\
 & \left.\times\dbinom{\lambda}{r_{1}+1}\cdots\dbinom{\lambda}{r_{i}+1}S_{1}\left(l+m+j,m+j\right)\right\} x^{j}.\nonumber
\end{align}

From (\ref{eq:28}), we note that
\begin{equation}
\left(1+e^{-\lambda t}\right)^{\mu}\frac{1}{\lif[k]\left(-t\right)}CP_{n}^{\left(k\right)}\left(x;\lambda,\mu\right)\sim\left(1,e^{-t}-1\right)\label{eq:39}
\end{equation}
and
\begin{equation}
x^{n}\sim\left(1,t\right).\label{eq:40}
\end{equation}

By (\ref{eq:24}), (\ref{eq:39}) and (\ref{eq:40}), we get
\begin{align}
 & \left(1+e^{-\lambda t}\right)^{\mu}\frac{1}{\lif[k]\left(-t\right)}CP_{n}^{\left(k\right)}\left(x;\lambda,\mu\right)\label{eq:41}\\
= & x\left(\frac{t}{e^{-t}-1}\right)^{n}x^{n-1}\nonumber \\
= & \left(-1\right)^{n}x\left(\frac{-t}{e^{-t}-1}\right)^{n}x^{n-1}\nonumber \\
= & \left(-1\right)^{n}\sum_{l=0}^{n-1}\left(-1\right)^{l}B_{l}^{\left(n\right)}\dbinom{n-1}{l}x^{n-l}.\nonumber
\end{align}

Thus, by (\ref{eq:41}), we see that
\begin{align}
 & CP_{n}^{\left(k\right)}\left(x;\lambda,\mu\right)\label{eq:42}\\
= & \left(-1\right)^{n}\sum_{l=0}^{n-1}\left(-1\right)^{l}\dbinom{n-1}{l}B_{l}^{\left(n\right)}\left(1+e^{-\lambda t}\right)^{-\mu}\lif[k]\left(-t\right)x^{n-l}\nonumber \\
= & \left(-1\right)^{n}\sum_{l=0}^{n-1}\left(-1\right)^{l}\dbinom{n-1}{l}B_{l}^{\left(n\right)}\sum_{m=0}^{n-l}\frac{\left(-1\right)^{m}\tbinom{n-l}{m}}{\left(m+1\right)^{k}}\left(1+e^{-\lambda t}\right)^{-\mu}x^{n-l-m}\nonumber \\
= & \left(-1\right)^{n}\sum_{l=0}^{n-1}\left(-1\right)^{l}\dbinom{n-1}{l}B_{l}^{\left(n\right)}\sum_{m=0}^{n-l}\frac{\left(-1\right)^{m}\tbinom{n-l}{m}}{\left(m+1\right)^{k}}\sum_{i=0}^{\infty}\sum_{j=0}^{\infty}\sum_{j_{1}+\cdots+j_{n}=j}2^{-\mu-i}\dbinom{-\mu}{i}\nonumber \\
 & \times\dbinom{j+i}{j_{1}+1,\cdots,j_{i}+1}\frac{\left(-\lambda t\right)^{j+i}}{\left(j+i\right)!}x^{n-l-m}\nonumber \\
= & \left(-1\right)^{n}\sum_{l=0}^{n}\sum_{m=0}^{n-l}\sum_{i=0}^{n-l-m}\sum_{j=0}^{n-l-m-i}\sum_{j_{1}+\cdots+j_{i}=n-l-m-i-r}\left(-1\right)^{n-r}\frac{2^{-\mu-i}\lambda^{n-l-m-r}}{\left(m+1\right)^{k}}\dbinom{n-1}{l}\nonumber \\
 & \times\dbinom{n-l}{m}\dbinom{-\mu}{i}\dbinom{n-l-m-r}{j_{1}+1,\cdots,j_{i}+1}\dbinom{n-l-m}{r}B_{l}^{\left(n\right)}x^{r}.\nonumber
\end{align}

Therefore, by (\ref{eq:42}), we obtain the following theorem.\end{rem*}
\begin{thm}
\label{thm:4}For $n\ge0$, we have
\begin{align*}
 & CP_{n}^{\left(k\right)}\left(x;\lambda,\mu\right)\\
= & \frac{\lambda^{n}}{2^{\mu}}\sum_{r=0}^{n}\left(-\lambda^{-1}\right)^{r}\left\{ \sum_{l=0}^{n-r}\sum_{m=0}^{n-r-l}\sum_{i=0}^{n-r-l-m}\sum_{j_{1}+\cdots+j_{i}=n-r-l-m-i}\frac{2^{-i}\lambda^{-l-m}}{\left(m+1\right)^{k}}\dbinom{n-1}{l}\right.\\
 & \left.\times\dbinom{n-l}{m}\dbinom{-\mu}{i}\dbinom{n-r-l-m}{j_{1}+1,\cdots,j_{i}+1}\dbinom{n-l-m}{r}B_{l}^{\left(n\right)}\right\} x^{r}.
\end{align*}

\end{thm}
$\,$
\begin{rem*}
By the same method as Theorem \ref{thm:4}, we get
\begin{align}
 & \hat{C}P_{n}^{\left(k\right)}\left(x;\lambda,\mu\right)\label{eq:43}\\
= & \frac{\lambda^{n}}{2^{\mu}}\sum_{r=0}^{n}\lambda^{-r}\left\{ \sum_{l=0}^{n-r}\sum_{m=0}^{n-l-r}\sum_{i=0}^{n-l-m-r}\sum_{j_{1}+\cdots+j_{i}=n-r-l-m-i}\frac{\left(-1\right)^{m}2^{-i}\lambda^{-l-m}}{\left(m+1\right)^{k}}\dbinom{n-1}{l}\right.\nonumber \\
 & \left.\times\dbinom{n-l}{m}\dbinom{-\mu}{i}\dbinom{n-r-l-m}{j_{1}+1,\cdots,j_{i}+1}\dbinom{n-l-m}{r}B_{l}^{\left(m\right)}\right\} x^{r}.\nonumber
\end{align}

\end{rem*}
$\,$

From (\ref{eq:12}), we note that
\begin{equation}
x^{\left(n\right)}=x\left(x+1\right)\cdots\left(x+n-1\right)\sim\left(1,1-e^{-t}\right).\label{eq:44}
\end{equation}

Thus, by (\ref{eq:44}), we see that
\begin{equation}
\left(-1\right)^{n}x^{\left(n\right)}=\left(-x\right)_{n}=\sum_{m=0}^{n}S_{1}\left(n,m\right)\left(-x\right)^{m}\sim\left(1,e^{-t}-1\right),\label{eq:45}
\end{equation}
and
\begin{equation}
\left(1+e^{-\lambda t}\right)^{\mu}\frac{1}{\lif[k]\left(-t\right)}CP_{n}^{\left(k\right)}\left(x;\lambda,\mu\right)\sim\left(1,e^{-t}-1\right).\label{eq:46}
\end{equation}

From (\ref{eq:45}), and (\ref{eq:46}), we have
\begin{align}
 & \left(1+e^{-\lambda t}\right)^{\mu}\frac{1}{\lif[k]\left(-t\right)}CP_{n}^{\left(k\right)}\left(x;\lambda,\mu\right)\label{eq:47}\\
= & \left(-1\right)^{n}x^{\left(n\right)}\nonumber \\
= & \sum_{l=0}^{n}S_{1}\left(n,l\right)\left(-x\right)^{l}\nonumber
\end{align}

Thus, by (\ref{eq:47}), we get
\begin{align}
 & CP_{n}^{\left(k\right)}\left(x;\lambda,\mu\right)\label{eq:48}\\
= & \sum_{l=0}^{n}S_{1}\left(n,l\right)\left(-1\right)^{l}\left(1+e^{-\lambda t}\right)^{-\mu}\lif[k]\left(-t\right)x^{l}\nonumber \\
= & \sum_{l=0}^{n}S_{1}\left(n,l\right)\left(-1\right)^{l}\sum_{m=0}^{l}\frac{\left(-1\right)^{m}\tbinom{l}{m}}{\left(m+1\right)^{k}}\left(1+e^{-\lambda t}\right)^{-\mu}x^{l-m}\nonumber \\
= & \sum_{l=0}^{n}S_{1}\left(n,l\right)\left(-1\right)^{l}\sum_{m=0}^{l}\frac{\left(-1\right)^{m}\tbinom{l}{m}}{\left(m+1\right)^{k}}\sum_{i=0}^{\infty}\sum_{j=0}^{\infty}\sum_{j_{1}+\cdots+j_{i}=j}2^{-\mu-i}\nonumber \\
 & \times\dbinom{-\mu}{i}\dbinom{j+i}{j_{1}+1,\cdots,j_{i}+1}\frac{\left(-\lambda\right)^{j+i}}{\left(j+i\right)!}t^{j+i}x^{l-m}\nonumber \\
= & \sum_{l=0}^{n}\sum_{m=0}^{l}\sum_{i=0}^{l-m}\sum_{r=0}^{l-m-i}\sum_{j_{1}+\cdots+j_{i}=l-m-i-r}\left(-1\right)^{r}\frac{2^{-\mu-i}\lambda^{l-m-r}}{\left(m+1\right)^{k}}\nonumber \\
 & \times\dbinom{l}{m}\dbinom{-\mu}{i}\dbinom{l-m-r}{j_{1}+1,\cdots,j_{i}+1}\dbinom{l-m}{r}S_{1}\left(n,l\right)x^{r}\nonumber \\
= & 2^{-\mu}\sum_{r=0}^{n}\left(-\lambda^{-1}\right)^{r}\left\{ \sum_{l=r}^{n}\sum_{m=0}^{l-r}\sum_{i=0}^{l-r-m}\sum_{j_{1}+\cdots+j_{i}=l-r-m-i}\frac{2^{-i}\lambda^{l-m}}{\left(m+1\right)^{k}}\right.\nonumber \\
 & \left.\times\dbinom{l}{m}\dbinom{-\mu}{i}\dbinom{l-m-r}{j_{1}+1,\cdots,j_{i}+1}\dbinom{l-m}{r}S_{1}\left(n,l\right)\right\} x^{r}\nonumber
\end{align}

Therefore, by (\ref{eq:48}), we obtain the following theorem.
\begin{thm}
\label{thm:5}For $n\ge0$, we have
\begin{align*}
 & CP_{n}^{\left(k\right)}\left(x;\lambda,\mu\right)\\
= & 2^{-\mu}\sum_{r=0}^{n}\left(-\lambda^{-1}\right)^{r}\left\{ \sum_{l=r}^{n}\sum_{m=0}^{l-r}\sum_{i=0}^{l-r-m}\sum_{j_{1}+\cdots+j_{i}=l-r-m-i}\frac{2^{-i}\lambda^{l-m}}{\left(m+1\right)^{k}}\right.\\
 & \left.\times\dbinom{l}{m}\dbinom{-\mu}{i}\dbinom{l-m-r}{j_{1}+1,\cdots,j_{i}+1}\dbinom{l-m}{r}S_{1}\left(n,l\right)\right\} x^{r}.
\end{align*}

\end{thm}
$\,$

It is easy to see that
\begin{equation}
\left(1+e^{\lambda t}\right)^{\mu}\frac{1}{\lif[k]\left(-t\right)}\hat{C}P_{n}^{\left(k\right)}\left(x;\lambda,\mu\right)\sim\left(1,e^{t}-1\right)\label{eq:49}
\end{equation}
and
\begin{equation}
\left(x\right)_{n}=x\left(x-1\right)\cdots\left(x-n+1\right)=\sum_{l=0}^{n}S_{1}\left(n,l\right)x^{l}\sim\left(1,e^{t}-1\right).\label{eq:50}
\end{equation}

By the same method as Theorem \ref{thm:5}, we get
\begin{align}
 & \hat{C}P_{n}^{\left(k\right)}\left(x;\lambda,\mu\right)\label{eq:51}\\
= & 2^{-\mu}\sum_{r=0}^{n}\lambda^{-r}\left\{ \sum_{l=r}^{n}\sum_{m=0}^{l-r}\sum_{i=0}^{l-r-m}\sum_{j_{1}+\cdots+j_{i}=l-r-m-i}\frac{\left(-1\right)^{m}2^{-i}\lambda^{l-m}}{\left(m+1\right)^{k}}\right.\nonumber \\
 & \times\left.\dbinom{l}{m}\dbinom{-\mu}{i}\dbinom{l-m-r}{j_{1}+1,\cdots,j_{i}+1}\dbinom{l-m}{r}S_{1}\left(n,l\right)\right\} x^{r}.\nonumber
\end{align}

From (\ref{eq:20}) and (\ref{eq:28}), we have
\begin{align}
 & CP_{n}^{\left(k\right)}\left(x;\lambda,\mu\right)\label{eq:52}\\
= & \sum_{j=0}^{n}\frac{1}{j!}\left\langle \left.\left(1+\left(1+t\right)^{\lambda}\right)^{-\mu}\lif[k]\left(\log\left(1+t\right)\right)\left(-\log\left(1+t\right)\right)^{j}\right|x^{n}\right\rangle x^{j}.\nonumber
\end{align}

Now, we observe that
\begin{align}
 & \left\langle \left.\left(1+\left(1+t\right)^{\lambda}\right)^{-\mu}\lif[k]\left(\log\left(1+t\right)\right)\left(-\log\left(1+t\right)\right)^{j}\right|x^{n}\right\rangle \label{eq:53}\\
= & \left(-1\right)^{j}\left\langle \log\left(1+t\right)^{j}\left|\sum_{m=0}^{\infty}CP_{m}^{\left(k\right)}\left(0;\lambda,\mu\right)\frac{t^{m}}{m!}x^{n}\right.\right\rangle \nonumber \\
= & \left(-1\right)^{j}\sum_{m=0}^{n}\dbinom{n}{m}CP_{m}^{\left(k\right)}\left(0;\lambda,\mu\right)\left\langle \left.\left(\log\left(1+t\right)\right)^{j}\right|x^{n-m}\right\rangle \nonumber \\
= & \left(-1\right)^{j}\sum_{m=0}^{n}\dbinom{n}{m}CP_{m}^{\left(k\right)}\left(0;\lambda,\mu\right)j!S_{1}\left(n-m,j\right).\nonumber
\end{align}

Therefore, by (\ref{eq:52}) and (\ref{eq:53}), we obtain the following
theorem.
\begin{thm}
\label{thm:6}For $n\ge0$, we have
\[
CP_{n}^{\left(k\right)}\left(x;\lambda,\mu\right)=\sum_{j=0}^{n}\left(-1\right)^{j}\left\{ \sum_{m=0}^{n}\dbinom{n}{m}S_{1}\left(n-m,j\right)CP_{m}^{\left(k\right)}\left(0;\lambda,\mu\right)\right\} x^{j}.
\]
\end{thm}
\begin{rem*}
By the same method as Theorem \ref{thm:6}, we get
\begin{equation}
\hat{C}P_{n}^{\left(k\right)}\left(x;\lambda,\mu\right)=\sum_{j=0}^{n}\left\{ \sum_{m=0}^{n}\dbinom{n}{m}S_{1}\left(n-m,j\right)\hat{C}P_{n}^{\left(k\right)}\left(0;\lambda,\mu\right)\right\} x^{j}.\label{eq:54}
\end{equation}

From (\ref{eq:21}), we have
\begin{equation}
CP_{n}^{\left(k\right)}\left(x+y;\lambda,\mu\right)=\sum_{j=0}^{n}\left(-1\right)^{j}\dbinom{n}{j}CP_{n-j}^{\left(k\right)}\left(x;\lambda,\mu\right)y^{\left(j\right)}\label{eq:55}
\end{equation}
and
\begin{equation}
\hat{C}P_{n}^{\left(k\right)}\left(x+y;\lambda,\mu\right)=\sum_{j=0}^{n}\dbinom{n}{j}\hat{C}P_{n-j}^{\left(k\right)}\left(x;\lambda,\mu\right)\left(y\right)_{j}.\label{eq:56}
\end{equation}

\end{rem*}
$\,$

By (\ref{eq:22}) and (\ref{eq:28}), we get
\begin{equation}
\left(e^{-t}-1\right)CP_{n}^{\left(k\right)}\left(x;\lambda,\mu\right)=nCP_{n-1}^{\left(k\right)}\left(x;\lambda,\mu\right)\label{eq:57}
\end{equation}
and
\begin{equation}
\left(e^{-t}-1\right)CP_{n}^{\left(k\right)}\left(x;\lambda,\mu\right)=CP_{n}^{\left(k\right)}\left(x-1;\lambda,\mu\right)-CP_{n}^{\left(k\right)}\left(x;\lambda,\mu\right).\label{eq:58}
\end{equation}

Therefore, by (\ref{eq:57}) and (\ref{eq:58}), we obtain the following
theorem.
\begin{thm}
\label{thm:7}For $n\ge0$, we have
\[
CP_{n}^{\left(k\right)}\left(x-1;\lambda,\mu\right)-CP_{n}^{\left(k\right)}\left(x;\lambda,\mu\right)=nCP_{n-1}^{\left(k\right)}\left(x;\lambda,\mu\right).
\]
\end{thm}
\begin{rem*}
By the same method as Theorem \ref{thm:7}, we get
\begin{equation}
\hat{C}P_{n}^{\left(k\right)}\left(x+1;\lambda,\mu\right)-\hat{C}P_{n}^{\left(k\right)}\left(x;\lambda,\mu\right)=n\hat{C}P_{n-1}^{\left(k\right)}\left(x;\lambda,\mu\right).\label{eq:59}
\end{equation}

From (\ref{eq:22}), (\ref{eq:28}), and (\ref{eq:29}), we have
\begin{align}
& CP_{n+1}^{\left(k\right)}\left(x;1,\mu\right)\label{eq:60}\\
& =-xCP_{n}^{\left(k\right)}\left(x+1;1,\mu\right)+\mu\sum_{m=0}^{n}\left(-\frac{1}{2}\right)^{m+1}\left(n\right)_{m}CP_{n-m}^{\left(k\right)}\left(x;1,\mu\right)\nonumber \\
& +2^{-\mu}\sum_{r=0}^{n}\left(-1\right)^{r}\left\{ \sum_{m=r}^{n}\sum_{l=r}^{m}\sum_{i=0}^{l-r}\sum_{j_{1}+\cdots+j_{i}=l-i-r}\frac{2^{-i}}{\left(m-l+2\right)^{k}}\dbinom{m}{l}\right.\nonumber \\
& \left.\times\dbinom{-\mu}{i}\dbinom{l-r}{j_{1}+1,\cdots,j_{i}+1}\dbinom{l}{r}S_{1}\left(n,m\right)\right\} \left(x+1\right)^{r},\nonumber
\end{align}
and
\begin{align}
& \hat{C}P_{n+1}^{\left(k\right)}\left(x;1,\mu\right)\label{eq:61}\\
= & x\hat{C}P_{n}^{\left(k\right)}\left(x-1;1,\mu\right)+\mu\sum_{m=0}^{n}\left(-\frac{1}{2}\right)^{m+1}\left(n\right)_{m}\hat{C}P_{n-m}^{\left(k\right)}\left(x;1,\mu\right)\nonumber \\
& -2^{-\mu}\sum_{r=0}^{n}\left\{ \sum_{m=r}^{n}\sum_{l=r}^{m}\sum_{i=0}^{l-r}\sum_{j_{1}+\cdots+j_{i}=l-i-r}\frac{(-1)^{m-l}2^{-i}}{\left(m-l+2\right)^{k}}\dbinom{m}{l}\dbinom{-\mu}{i}\right.\nonumber \\
& \left.\times\dbinom{l-r}{j_{1}+1,\cdots,j_{i}+1}\dbinom{l}{r}S_{1}\left(n,m\right)\right\} \left(x-1\right)^{r}.\nonumber
\end{align}
\end{rem*}
$\,$
By (\ref{eq:14}) and (\ref{eq:27}), we get
\begin{align}
 & CP_{n}^{\left(k\right)}\left(y;\lambda,\mu\right)\label{eq:62}\\
= & \left\langle \left.\sum_{l=0}^{\infty}CP_{l}^{\left(k\right)}\left(y;\lambda,\mu\right)\frac{t^{l}}{l!}\right|x^{n}\right\rangle \nonumber \\
= & \left\langle \left.\left(1+\left(1+t\right)^{\lambda}\right)^{-\mu}\lif[k]\left(\log\left(1+t\right)\right)\left(1+t\right)^{-y}\right|x\cdot x^{n-1}\right\rangle \nonumber \\
= & \left\langle \left.\partial_{t}\left(\left(1+\left(1+t\right)^{\lambda}\right)^{-\mu}\lif[k]\left(\log\left(1+t\right)\right)\left(1+t\right)^{-y}\right)\right|x^{n-1}\right\rangle \nonumber \\
= & \left\langle \left.\left(\partial_{t}\left(1+\left(1+t\right)^{\lambda}\right)^{-\mu}\right)\lif[k]\left(\log\left(1+t\right)\right)\left(1+t\right)^{-y}\right|x^{n-1}\right\rangle \nonumber \\
 & +\left\langle \left.\left(1+\left(1+t\right)^{\lambda}\right)^{-\mu}\left(\partial_{t}\lif[k]\left(\log\left(1+t\right)\right)\right)\left(1+t\right)^{-y}\right|x^{n-1}\right\rangle \nonumber \\
 & +\left\langle \left.\left(1+\left(1+t\right)^{\lambda}\right)^{-\mu}\lif[k]\left(\log\left(1+t\right)\right)\left(\partial_{t}\left(1+t\right)^{-y}\right)\right|x^{n-1}\right\rangle \nonumber \\
= & -\mu\lambda\left\langle \left.\left(1+\left(1+t\right)^{\lambda}\right)^{-\mu-1}\lif[k]\left(\log\left(1+t\right)\right)\left(1+t\right)^{-\left(y-\lambda+1\right)}\right|x^{n-1}\right\rangle \nonumber \\
 & -y\left\langle \left.\left(1+\left(1+t\right)^{\lambda}\right)^{-\mu}\left(\lif[k]\left(\log\left(1+t\right)\right)\right)\left(1+t\right)^{-y-1}\right|x^{n-1}\right\rangle \nonumber \\
 & +\left\langle \left.\left(1+\left(1+t\right)^{\lambda}\right)^{-\mu}\left(\partial_{t}\lif[k]\left(\log\left(1+t\right)\right)\right)\left(1+t\right)^{-y}\right|x^{n-1}\right\rangle \nonumber \\
= & -\mu\lambda CP_{n-1}^{\left(k\right)}\left(y-\lambda+1;\lambda,\mu+1\right)-yCP_{n-1}^{\left(k\right)}\left(y+1;\lambda,\mu\right)\nonumber \\
 & +\left\langle \left.\left(1+\left(1+t\right)^{\lambda}\right)^{-\mu}\frac{\lif[k-1]\left(\log\left(1+t\right)\right)-\lif[k]\left(\log\left(1+t\right)\right)}{\left(1+t\right)\log\left(1+t\right)}\left(1+t\right)^{-y}\right|x^{n-1}\right\rangle .\nonumber
\end{align}

Now, we observe that
\begin{align}
 & \left\langle \left.\left(1+\left(1+t\right)^{\lambda}\right)^{-\mu}\frac{\lif[k-1]\left(\log\left(1+t\right)\right)-\lif[k]\left(\log\left(1+t\right)\right)}{\left(1+t\right)\log\left(1+t\right)}\left(1+t\right)^{-y}\right|x^{n-1}\right\rangle \label{eq:63}\\
= & \left\langle \left(1+\left(1+t\right)^{\lambda}\right)^{-\mu}\frac{\lif[k-1]\left(\log\left(1+t\right)\right)-\lif[k]\left(\log\left(1+t\right)\right)}{t}\left(1+t\right)^{-y}\right|\nonumber \\
 & \left.\frac{t}{\left(1+t\right)\log\left(1+t\right)}x^{n-1}\right\rangle \nonumber \\
= & \sum_{l=0}^{n-1}\dbinom{n-1}{l}\hat{\mathbb{C}}_{n-1-l}^{\left(1\right)}\left(0\right)\nonumber \\
 & \times\left\langle \left.\left(1+\left(1+t\right)^{\lambda}\right)^{-\mu}\frac{\lif[k-1]\left(\log\left(1+t\right)\right)-\lif[k]\left(\log\left(1+t\right)\right)}{t}\left(1+t\right)^{-y}\right|x^{l}\right\rangle \nonumber \\
= & \sum_{l=0}^{n-1}\dbinom{n-1}{l}\hat{\mathbb{C}}_{n-1-l}^{\left(1\right)}\left(0\right)\nonumber \\
 & \times\left\langle \left.\left(1+\left(1+t\right)^{\lambda}\right)^{-\mu}\frac{\lif[k-1]\left(\log\left(1+t\right)\right)-\lif[k]\left(\log\left(1+t\right)\right)}{t}\left(1+t\right)^{-y}\right|t\frac{x^{l+1}}{l+1}\right\rangle \nonumber \\
= & \frac{1}{n}\sum_{l=0}^{n-1}\dbinom{n}{l+1}\hat{\mathbb{C}}_{n-1-l}^{\left(1\right)}\left(0\right)\left\{ CP_{l+1}^{\left(k-1\right)}\left(y;\lambda,\mu\right)-CP_{l+1}^{\left(k\right)}\left(y;\lambda,\mu\right)\right\} .\nonumber
\end{align}

Therefore, by (\ref{eq:62}) and (\ref{eq:63}), we obtain the following
theorem.
\begin{thm}
\label{thm:8}For $n\ge0$, we have
\begin{align*}
 & CP_{n}^{\left(k\right)}\left(x;\lambda,\mu\right)\\
= & -\mu\lambda CP_{n-1}^{\left(k\right)}\left(x-\lambda+1;\lambda,\mu+1\right)-xCP_{n-1}^{\left(k\right)}\left(x+1;\lambda,\mu\right)\\
 & +\frac{1}{n}\sum_{l=0}^{n-1}\dbinom{n}{l+1}\hat{C}_{n-1-l}\left\{ CP_{l+1}^{\left(k-1\right)}\left(x;\lambda,\mu\right)-CP_{l+1}^{\left(k\right)}\left(x;\lambda,\mu\right)\right\} .
\end{align*}
where $\hat{C}_{n-1-l}=\hat{\mathbb{C}}_{n-1-l}^{\left(1\right)}\left(0\right).$\end{thm}
\begin{rem*}
By the same method as Theorem \ref{thm:8}, we get
\begin{align}
 & \hat{C}P_{n}^{\left(k\right)}\left(x;\lambda,\mu\right)\label{eq:64}\\
= & -\mu\lambda\hat{C}P_{n-1}^{\left(k\right)}\left(x+\lambda-1;\lambda,\mu+1\right)+x\hat{C}P_{n-1}^{\left(k\right)}\left(x-1;\lambda,\mu\right)\nonumber \\
 & +\frac{1}{n}\sum_{l=0}^{n-1}\dbinom{n}{l+1}\hat{C}_{n-1-l}\left(\hat{C}P_{l+1}^{\left(k-1\right)}\left(x;\lambda,\mu\right)-\hat{C}P_{l+1}^{\left(k\right)}\left(x;\lambda,\mu\right)\right).\nonumber
\end{align}

\end{rem*}
$\,$

By (\ref{eq:23}), we get
\begin{align}
 & \frac{d}{dx}CP_{n}^{\left(k\right)}\left(x;\lambda,\mu\right)\label{eq:65}\\
= & \sum_{l=0}^{n-1}\dbinom{n}{l}\left\langle \left.-\log\left(1+t\right)\right|x^{n-l}\right\rangle CP_{l}^{\left(k\right)}\left(x;\lambda,\mu\right)\nonumber \\
= & \sum_{l=0}^{n-1}\dbinom{n}{l}\left\langle \left.\sum_{m=1}^{\infty}\frac{\left(-1\right)^{m}}{m}t^{m}\right|x^{n-l}\right\rangle CP_{l}^{\left(k\right)}\left(x;\lambda,\mu\right)\nonumber \\
= & \sum_{l=0}^{n-1}\dbinom{n}{l}\sum_{m=1}^{\infty}\frac{\left(-1\right)^{m}}{m}\left\langle \left.t^{m}\right|x^{n-l}\right\rangle CP_{l}^{\left(k\right)}\left(x;\lambda,\mu\right)\nonumber \\
= & \sum_{l=0}^{n-1}\dbinom{n}{l}\left(-1\right)^{n-l}CP_{l}^{\left(k\right)}\left(x;\lambda,\mu\right)\left(n-l-1\right)!\nonumber \\
= & n!\sum_{l=0}^{n-1}\frac{\left(-1\right)^{n-l}}{\left(n-l\right)l!}CP_{l}^{\left(k\right)}\left(x;\lambda,\mu\right).\nonumber
\end{align}

By the same method as (\ref{eq:65}), we get
\begin{align}
 & \frac{d}{dx}\hat{C}P_{n}^{\left(k\right)}\left(x;\lambda,\mu\right)\label{eq:66}\\
= & n!\sum_{l=0}^{n-1}\frac{\left(-1\right)^{n-l-1}}{\left(n-l\right)l!}\hat{C}P_{l}^{\left(k\right)}\left(x;\lambda,\mu\right).\nonumber
\end{align}

Now, we compute the following equation in two different ways:

\[
\left\langle \left.\left(1+\left(1+t\right)^{\lambda}\right)^{-\mu}\lif[k]\left(-\log\left(1+t\right)\right)\left(\log\left(1+t\right)\right)^{m}\right|x^{n}\right\rangle .
\]

On the one hand,

\begin{align}
 & \left\langle \left.\left(1+\left(1+t\right)^{\lambda}\right)^{-\mu}\lif[k]\left(-\log\left(1+t\right)\right)\left(\log\left(1+t\right)\right)^{m}\right|x^{n}\right\rangle \label{eq:67}\\
= & \left\langle \left.\left(1+\left(1+t\right)^{\lambda}\right)^{-\mu}\lif[k]\left(-\log\left(1+t\right)\right)\right|\left(\log\left(1+t\right)\right)^{m}x^{n}\right\rangle \nonumber \\
= & \sum_{l=0}^{n-m}m!\dbinom{n}{l+m}S_{1}\left(l+m,m\right)\left\langle \left.\left(1+\left(1+t\right)^{\lambda}\right)^{-\mu}\lif[k]\left(-\log\left(1+t\right)\right)\right|x^{n-l-m}\right\rangle \nonumber \\
= & \sum_{l=0}^{n-m}m!\dbinom{n}{l}S_{1}\left(n-l,m\right)\hat{C}P_{l}^{\left(k\right)}\left(0;\lambda,\mu\right).\nonumber
\end{align}

On the other hand,
\begin{align}
 & \left\langle \left.\left(1+\left(1+t\right)^{\lambda}\right)^{-\mu}\lif[k]\left(-\log\left(1+t\right)\right)\left(\log\left(1+t\right)\right)^{m}\right|x^{n}\right\rangle \label{eq:68}\\
= & \left\langle \left.\partial_{t}\left(\left(1+\left(1+t\right)^{\lambda}\right)^{-\mu}\lif[k]\left(-\log\left(1+t\right)\right)\left(\log\left(1+t\right)\right)^{m}\right)\right|x^{n-1}\right\rangle \nonumber \\
= & \left\langle \left.\left(\partial_{t}\left(1+\left(1+t\right)^{\lambda}\right)^{-\mu}\right)\lif[k]\left(-\log\left(1+t\right)\right)\left(\log\left(1+t\right)\right)^{m}\right|x^{n-1}\right\rangle \nonumber \\
 & +\left\langle \left.\left(1+\left(1+t\right)^{\lambda}\right)^{-\mu}\left(\partial_{t}\lif[k]\left(-\log\left(1+t\right)\right)\right)\left(\log\left(1+t\right)\right)^{m}\right|x^{n-1}\right\rangle \nonumber \\
 & +\left\langle \left.\left(1+\left(1+t\right)^{\lambda}\right)^{-\mu}\lif[k]\left(-\log\left(1+t\right)\right)\left(\partial_{t}\left(\log\left(1+t\right)\right)^{m}\right)\right|x^{n-1}\right\rangle .\nonumber
\end{align}

Therefore, by (\ref{eq:67}) and (\ref{eq:68}), we obtain the following
theorem.
\begin{thm}
\label{thm:9}For $n\in\mathbb{N}$ with $n\ge2$, let $n-1\ge m\ge1.$
Then we have
\begin{align*}
 & m\sum_{l=0}^{n-m}\dbinom{n}{l}S_{1}\left(n-l,m\right)\hat{C}P_{l}^{\left(k\right)}\left(0;\lambda,\mu\right)\\
= & -\mu\lambda m\sum_{l=0}^{n-1-m}\dbinom{n-1}{l}S_{1}\left(n-1-l,m\right)\hat{C}P_{l}^{\left(k\right)}\left(\lambda-1;\lambda,\mu+1\right)\\
 & +\sum_{l=0}^{n-m}\dbinom{n-1}{l}S_{1}\left(n-1-l,m-1\right)\hat{C}P_{l}^{\left(k-1\right)}\left(-1;\lambda,\mu\right)\\
 & +\left(m-1\right)\sum_{l=0}^{n-m}\dbinom{n-1}{l}S_{1}\left(n-1-l,m-1\right)\hat{C}P_{l}^{\left(k\right)}\left(-1;\lambda,\mu\right).
\end{align*}
\end{thm}
\begin{rem*}
By the same method as Theorem \ref{thm:9}, we get
\begin{align*}
 & m\sum_{l=0}^{n-m}\dbinom{n}{l}S_{1}\left(n-l,m\right)CP_{l}^{\left(k\right)}\left(0;\lambda,\mu\right)\\
= & -\mu\lambda m\sum_{l=0}^{n-1-m}\dbinom{n-1}{l}S_{1}\left(n-1-l,m\right)CP_{l}^{\left(k\right)}\left(1-\lambda;\lambda,\mu+1\right)\\
 & +\sum_{l=0}^{n-m}\dbinom{n-1}{l}S_{1}\left(n-1-l,m-1\right)CP_{l}^{\left(k-1\right)}\left(1;\lambda,\mu\right)\\
 & +\left(m-1\right)\sum_{l=0}^{n-m}\dbinom{n-1}{l}S_{1}\left(n-1-l,m-1\right)CP_{l}^{\left(k\right)}\left(-1;\lambda,\mu\right),
\end{align*}
where $n-1\ge m\ge1$.
\end{rem*}
$\,$

Let us consider the following two sheffer sequences :

\begin{equation}
CP_{n}^{\left(k\right)}\left(x;\lambda,\mu\right)\sim\left(\left(1+e^{-\lambda t}\right)^{\mu}\frac{1}{\lif[k]\left(-t\right)},e^{-t}-1\right)\label{eq:69}
\end{equation}
and
\begin{equation}
B_{n}^{\left(s\right)}\left(x\right)\sim\left(\left(\frac{e^{t}-1}{t}\right)^{s},t\right),\quad\left(s\in\mathbb{Z}_{\ge0}\right).\label{eq:70}
\end{equation}

Let
\begin{equation}
CP_{n}^{\left(k\right)}\left(x;\lambda,\mu\right)=\sum_{m=0}^{n}C_{n,m}B_{m}^{\left(s\right)}\left(x\right).\label{eq:71}
\end{equation}

Then, by (\ref{eq:26}), we get
\begin{align}
 & C_{n,m}\label{eq:72}\\
= &\frac{1}{m!} \left\langle \left.\frac{\left(\frac{e^{-\log\left(1+t\right)}-1}{-\log\left(1+t\right)}\right)^{s}}{\left(1+e^{\lambda\log\left(1+t\right)}\right)^{\mu}}\lif[k]\left(\log\left(1+t\right)\right)\left(-\log\left(1+t\right)\right)^{m}\right|x^{n}\right\rangle \nonumber \\
= & \frac{\left(-1\right)^{m}}{m!}\left\langle \left.\left(1+\left(1+t\right)^{\lambda}\right)^{-\mu}\lif[k]\left(\log\left(1+t\right)\right)\right.\right. \nonumber \\
 & \times\left.\left.\left(1+t\right)^{-s}\left(\frac{t}{\log\left(1+t\right)}\right)^{s}\right|\left(\log\left(1+t\right)\right)^{m}x^{n}\right\rangle \nonumber \\
= & \frac{\left(-1\right)^{m}}{m!}\sum_{l=0}^{n-m}m!\dbinom{n}{l+m}S_{1}\left(l+m,m\right)\sum_{i=0}^{n-l-m}\dbinom{n-l-m}{i}\mathbb{C}_{i}^{\left(s\right)}\nonumber \\
 & \times\left\langle \left.\left(1+\left(1+t\right)^{\lambda}\right)^{-\mu}\lif[k]\left(\log\left(1+t\right)\right)\left(1+t\right)^{-s}\right|x^{n-l-m-i}\right\rangle \nonumber \\
= & \left(-1\right)^{m}\sum_{l=0}^{n-m}\dbinom{n}{l}S_{1}\left(n-l,m\right)\sum_{i=0}^{l}\dbinom{l}{i}\mathbb{C}_{i}^{\left(s\right)}CP_{l-i}^{\left(k\right)}\left(s;\lambda,\mu\right).\nonumber
\end{align}

Therefore, by (\ref{eq:71}) and (\ref{eq:72}), we obtain the following
the theorem.
\begin{thm}
\label{thm:10}For $n\ge0$, we have
\begin{align*}
 & CP_{n}^{\left(k\right)}\left(x;\lambda,\mu\right)\\
= & \sum_{m=0}^{n}\left(-1\right)^{m}\left\{ \sum_{l=0}^{n-m}\sum_{i=0}^{l}\dbinom{n}{l}\dbinom{l}{i}S_{1}\left(n-l,m\right)\mathbb{C}_{i}^{\left(s\right)}CP_{l-i}^{\left(k\right)}\left(s;\lambda,\mu\right)\right\} B_{n}^{\left(s\right)}\left(x\right).
\end{align*}
\end{thm}
\begin{rem*}
By the same method as Theorem \ref{thm:10}, we have
\begin{align}
 & \hat{C}P_{n}^{\left(k\right)}\left(x;\lambda,\mu\right)\label{eq:73}\\
= & \sum_{m=0}^{n}\left\{ \sum_{l=0}^{n-m}\sum_{i=0}^{l}\dbinom{n}{l}\dbinom{l}{i}S_{1}\left(n-l,m\right)\mathbb{\hat{C}}_{i}^{\left(s\right)}\hat{C}P_{l-i}^{\left(k\right)}\left(s;\lambda,\mu\right)\right\} B_{m}^{\left(s\right)}\left(x\right).\nonumber
\end{align}

\end{rem*}
$\,$

For $CP_{n}^{\left(k\right)}\left(x;\lambda,\mu\right)\sim\left(\left(1+e^{-\lambda t}\right)^{\mu}\frac{1}{\lif[k]\left(-t\right)},e^{-t}-1\right)$,
$H_{n}^{\left(s\right)}\left(x|\lambda\right)\sim\left(\left(\frac{e^{t}-\lambda}{1-\lambda}\right)^{s},t\right)$,
$s\in\mathbb{Z}_{\ge0}$, $\lambda\in\mathbb{C}$ with $\lambda\ne1$,
let us assume that
\begin{equation}
CP_{n}^{\left(k\right)}\left(x;\lambda,\mu\right)=\sum_{m=0}^{n}C_{n,m}H_{m}^{\left(s\right)}\left(x;\lambda\right).\label{eq:74}
\end{equation}

From (\ref{eq:26}), we have
\begin{align}
 & C_{n,m}\label{eq:75}\\
= & \frac{\left(-1\right)^{m}}{m!}\left\langle \left(1+\left(1+t\right)^{\lambda}\right)^{-\mu}\lif[k]\left(\log\left(1+t\right)\right)\right.\nonumber \\
 & \times\left.\left.\left(1+t\right)^{-s}\left(1+\frac{\lambda}{\lambda-1}t\right)^{s}\right|\left(\log\left(1+t\right)\right)^{m}x^{n}\right\rangle \nonumber \\
= & \frac{\left(-1\right)^{m}}{m!}\sum_{l=0}^{n-m}m!\dbinom{n}{l+m}S_{1}\left(l+m,m\right)\sum_{i=0}^{\min\left\{ s,n-l-m\right\} }\dbinom{s}{i}\left(\frac{\lambda}{\lambda-1}\right)^{i}\nonumber \\
 & \times\left\langle \left.\left(1+\left(1+t\right)^{\lambda}\right)^{-\mu}\lif[k]\left(\log\left(1+t\right)\right)\left(1+t\right)^{-s}\right|t^{i}x^{n-l-m}\right\rangle \nonumber \\
= & \left(-1\right)^{m}\sum_{l=0}^{n-m}\sum_{i=0}^{\min\left\{ s,n-l-m\right\} }\dbinom{n}{l+m}\dbinom{s}{i}\nonumber \\
 & \times\left(n-l-m\right)_{i}\left(\frac{\lambda}{\lambda-1}\right)^{i}S_{1}\left(l+m,m\right)CP_{n-l-m-i}^{\left(k\right)}\left(s;\lambda,\mu\right)\nonumber \\
= & \left(-1\right)^{m}\sum_{l=0}^{n-m}\sum_{i=0}^{\min\left\{ s,l\right\} }\dbinom{n}{l}\dbinom{s}{i}\left(l\right)_{i}\left(\frac{\lambda}{\lambda-1}\right)^{i}S_{1}\left(n-l,m\right)CP_{l-i}^{\left(k\right)}\left(s;\lambda,\mu\right).\nonumber
\end{align}

Therefore, by (\ref{eq:75}) and (\ref{eq:76}), we obtain the following
theorem.
\begin{thm}
\label{thm:11}For $\lambda\in\mathbb{C}$ with $\lambda\ne1$, $n\ge0$,
we have
\begin{align*}
 & CP_{n}^{\left(k\right)}\left(x;\lambda,\mu\right)\\
= & \sum_{m=0}^{n}\left(-1\right)^{m}\left\{ \sum_{l=0}^{n-m}\sum_{i=0}^{\min\left\{ s,l\right\} }\dbinom{n}{l}\dbinom{s}{i}\left(l\right)_{i}.\right.\\
 & \left.\left(\frac{\lambda}{\lambda-1}\right)^{i}S_{1}\left(n-l,m\right)CP_{l-i}^{\left(k\right)}\left(s;\lambda,\mu\right)\right\} H_{m}^{\left(s\right)}\left(x;\lambda\right).
\end{align*}
\end{thm}
\begin{rem*}
By the same method as Theorem \ref{thm:11}, we get

\begin{align*}
 & \hat{C}P_{n}^{\left(k\right)}\left(x;\lambda,\mu\right)\\
= & \sum_{m=0}^{n}\left\{ \sum_{l=0}^{n-m}\sum_{i=0}^{\min\left\{ s,l\right\} }\dbinom{n}{l}\dbinom{s}{i}\left(l\right)_{i}.\right.\\
 & \left.\left(\frac{1}{1-\lambda}\right)^{i}S_{1}\left(n-l,m\right)\hat{C}P_{l-i}^{\left(k\right)}\left(0;\lambda,\mu\right)\right\} H_{m}^{\left(s\right)}\left(x;\lambda\right)
\end{align*}

\end{rem*}
$\,$

For $CP_{n}^{\left(k\right)}\left(x;\lambda,\mu\right)\sim\left(\left(1+e^{-\lambda t}\right)^{\mu}\frac{1}{\lif[k]\left(-t\right)},e^{-t}-1\right)$
and $x^{\left(n\right)}\sim\left(1,1-e^{-t}\right)$, let us assume
that
\begin{equation}
CP_{n}^{\left(k\right)}\left(x;\lambda,\mu\right)=\sum_{m=0}^{\infty}C_{n,m}x^{\left(m\right)}.\label{eq:76}
\end{equation}

By (\ref{eq:26}), we get
\begin{align}
 & C_{n,m}\label{eq:77}\\
= & \frac{1}{m!}\left\langle \left.\frac{1}{\left(1+e^{\lambda\log\left(1+t\right)}\right)^{\mu}}\lif[k]\left(\log\left(1+t\right)\right)\left(1-e^{\log\left(1+t\right)}\right)^{m}\right|x^{n}\right\rangle \nonumber \\
= & \frac{1}{m!}\left\langle \left.\left(1+\left(1+t\right)^{\lambda}\right)^{-\mu}\lif[k]\left(\log\left(1+t\right)\right)\left(-t\right)^{m}\right|x^{n}\right\rangle \nonumber \\
= & \frac{\left(-1\right)^{m}}{m!}\left\langle \left.\left(1+\left(1+t\right)^{\lambda}\right)^{-\mu}\lif[k]\left(\log\left(1+t\right)\right)\right|t^{m}x^{n}\right\rangle \nonumber \\
= & \left(-1\right)^{m}\dbinom{n}{m}\left\langle \left.\left(1+\left(1+t\right)^{\lambda}\right)^{-\mu}\lif[k]\left(\log\left(1+t\right)\right)\right|x^{n-m}\right\rangle \nonumber \\
= & \left(-1\right)^{m}\dbinom{n}{m}CP_{n-m}^{\left(k\right)}\left(0;\lambda,\mu\right).\nonumber
\end{align}

Therefore, by (\ref{eq:77}) and (\ref{eq:78}), we obtain the following
theorem.
\begin{thm}
For $n\ge0$, we have
\[
CP_{n}^{\left(k\right)}\left(x;\lambda,\mu\right)=\sum_{m=0}^{n}\left(-1\right)^{m}\dbinom{n}{m}CP_{n-m}^{\left(k\right)}\left(0;\lambda,\mu\right)x^{\left(m\right)}.
\]
\end{thm}
\begin{rem*}
By the same method as Theorem 12, we get
\begin{equation}
\hat{C}P_{n}^{\left(k\right)}\left(x;\lambda,\mu\right)=\sum_{m=0}^{n}\dbinom{n}{m}\hat{C}P_{n-m}^{\left(k\right)}\left(0;\lambda,\mu\right)\left(x\right)_{m}.\label{eq:78}
\end{equation}
\end{rem*}

$\,$

\noindent \noun{Department of Mathematics, Sogang University, Seoul
121-742, Republic of Korea}

\noindent \emph{E-mail}\noun{ }\emph{address : }\texttt{dskim@sogang.ac.kr}

\noun{$\,$}

\noindent \noun{Department of Mathematics, Kwangwoon University, Seoul
139-701, Republic of Korea}

\noindent \emph{E-mail}\noun{ }\emph{address : }\texttt{tkkim@kw.ac.kr}
\end{document}